\documentstyle[12pt,newlfont]{amsart}

\newcommand{\nc}{\newcommand}

\nc{\one}{\mbox{\bf 1}}
\nc{\invtensor}{\underset{\leftarrow}{\otimes}}
\nc{\ad}{\operatorname{ad}}
\nc{\rk}{\operatorname{rank}}
\nc{\corank}{\operatorname{corank}}
\nc{\Sym}{\operatorname{Sym}}

\nc{\sym}{\operatorname{sym}}
\nc{\id}{\operatorname{id}}
\nc{\htt}{\operatorname{ht}}
\nc{\Ker}{\operatorname{Ker}}
\nc{\im}{\operatorname{Im}}
\nc{\re}{\operatorname{Re}}
\nc{\sn}{\operatorname{sn}}
\nc{\spn}{\operatorname{span}}

\nc{\Irr}{\operatorname{Irr}}
\nc{\sgn}{\operatorname{sgn}}
\nc{\F}{\operatorname{F}}

\nc{\Soc}{\operatorname{Soc}}
\nc{\Inj}{\operatorname{E}}
\nc{\Hom}{\operatorname{Hom}}
\nc{\End}{\operatorname{End}}
\nc{\supp}{\operatorname{supp}}
\nc{\Card}{\operatorname{Card}}
\nc{\Mod}{\operatorname{Mod}}
\nc{\Ann}{\operatorname{Ann}}
\nc{\Ind}{\operatorname{Ind}}
\nc{\Coind}{\operatorname{Coind}}

\nc{\wt}{\operatorname{wt}}
\nc{\ch}{\operatorname{ch}}
\nc{\Stab}{\operatorname{Stab}}
\nc{\Sch}{{\cal S}\mbox{\em ch}}

\nc{\Spec}{\operatorname{Spec}}
\nc{\Prim}{\operatorname{Prim}}
\nc{\Max}{\operatorname{Max}}
\nc{\Aut}{\operatorname{Aut}}
\nc{\Fract}{\operatorname{Fract}}
\nc{\gr}{\operatorname{gr}}
\nc{\Gr}{\operatorname{Gr}}
\nc{\wdM}{\widetilde{M}}
\nc{\wdV}{\widetilde{V}}
\nc{\wdN}{\widetilde{N}}
\nc{\wdchi}{\widetilde{\chi}}

\nc{\pprv}{{\cal J}}
\nc{\eprv}{{\det{\cal J}}}
\nc{\uu}{{\cal I}}

\nc{\cO}{\operatorname{\cal O}}
\nc{\wdO}{\operatorname{\widetilde{\cal O}}}
\nc{\Ob}{\operatorname{\cal Ob}}
\nc{\Dglie}{\operatorname{{\cal D}glie}}
\nc{\Fin}{\operatorname{{\cal F}in}}
\nc{\cC}{\operatorname{\cal C}}
\nc{\wdC}{\operatorname{\widetilde{\cal C}}}

\nc{\Sg}{{\cal S}({\frak g})}

\nc{\Ug}{\widetilde{\cal U}}
\nc{\Zg}{{\cal Z}({\frak g})}
\nc{\tZg}{{\widetilde{\cal Z}({\frak g})}}
\nc{\Zk}{{\cal Z}({\frak k})}
\nc{\Sh}{{\cal S}({\frak h})}
\nc{\Uh}{{\cal U}({\frak h})}

\nc{\Uk}{{\cal U}({\frak k})}
\nc{\Ag}{{\cal A}({\frak g})}

\nc{\cZ}{\cal Z}
\nc{\cS}{\cal S}
\nc{\cP}{\cal P}
\nc{\cL}{\cal L}
\nc{\cU}{\cal U}
\nc{\cH}{\cal H}
\nc{\cK}{\cal K}
\nc{\cF}{\cal F}
\nc{\fg}{\frak g}
\nc{\CO}{\cal O}
\nc{\fn}{\frak n}
\nc{\fm}{\frak m}
\nc{\fh}{\frak h}
\nc{\ft}{\frak t}
\nc{\fk}{\frak k}
\nc{\fp}{\frak p}
\nc{\fI}{\frak I}
\nc{\veps}{\varepsilon}
\nc{\fsl}{\frak {sl}}

\nc{\dirlim}{\underset{\rightarrow}{\lim}\,} 
\nc{\nen}{\newenvironment}
\nc{\ol}{\overline}
\nc{\ul}{\underline}
\nc{\ra}{\rightarrow}
\nc{\lra}{\longrightarrow}
\nc{\Lra}{\Longrightarrow}
\nc{\Lla}{\Longleftarrow}
\nc{\Llra}{\Longleftrightarrow}
\nc{\thla}{\twoheadleftarrow}
\nc{\hra}{\hookrightarrow}
\nc{\iso}{\overset{\sim}{\lra}}
\nc{\ssubset}{\underset{\not=}{\subset}}

\nc{\Thm}[1]{Theorem~\ref{#1}}
\nc{\Prop}[1]{Proposition~\ref{#1}}
\nc{\Lem}[1]{Lemma~\ref{#1}}
\nc{\Cor}[1]{Corollary~\ref{#1}}
\nc{\Conj}[1]{Conjecture~\ref{#1}}
\nc{\Claim}[1]{Claim~\ref{#1}}
\nc{\Defn}[1]{Definition~\ref{#1}}
\nc{\Exa}[1]{Example~\ref{#1}}
\nc{\Rem}[1]{Remark~\ref{#1}}
\nc{\Note}[1]{Note~\ref{#1}}
\nc{\Quest}[1]{Question~\ref{#1}}
\nc{\Hyp}[1]{Hypoth\`ese~\ref{#1}}

\nen{thm}[1]{\label{#1}{\bf Theorem.\ } \em}{}
\nen{prop}[1]{\label{#1}{\bf Proposition.\ } \em}{}
\nen{lem}[1]{\label{#1}{\bf Lemma.\ } \em}{}
\nen{cor}[1]{\label{#1}{\bf Corollary.\ } \em}{}
\nen{conj}[1]{\label{#1}{\bf Conjecture.\ } \em}{}
\nen{claim}[1]{\label{#1}{\bf Claim.\ } \em}{}

\nen{defn}[1]{\label{#1}{\bf Definition.\ } }{}
\nen{exa}[1]{\label{#1}{\bf Example.\ } }{}

\nen{rem}[1]{\label{#1}{\em Remark.\ } }{}
\nen{note}[1]{\label{#1}{\em Note.\ } }{}
\nen{exer}[1]{\label{#1}{\em Exercise.\ } }{}
\nen{sket}[1]{\label{#1}{\em Sketch of proof.\ } }{}
\nen{quest}[1]{\label{#1}{\bf Question.\ } \em}{}
\nen{hyp}[1]{\label{#1}{\bf Hypoth\`ese.\ } \em}{}
\begin{document}


\title[]{Strongly typical representations of the basic classical Lie
superalgebras.}

\author[]{Maria Gorelik}

\address{ 
{\tt email: gorelik@@mpim-bonn.mpg.de} 
}

\thanks{The author was partially supported by TMR Grant No. FMRX-CT97-0100.
Research at MSRI is supported in part by NSF grant DMS-9701755.}

\begin{abstract}
We describe the category of representations with a 
strongly typical central character of a basic classical Lie superalgebra
in terms of representations of its even part.
\end{abstract}

\maketitle

\section{Introduction}

\subsection{}
In~\cite{psI}, I.~Penkov and V.~Serganova show that
the category of representations of a  basic classical Lie superalgebra $\fg$
of type I with a fixed typical central character
is equivalent to the category of representations of the even part $\fg_0$
with a suitable central character. A similar result for type II
was proven by I.~Penkov in~\cite{p} for ``generic'' central
characters. The aim of this paper is to understand for which
central characters such an equivalence holds. We also
study a simplest example when such an equivalence fails to exist,
but the corresponding category of representations of a Lie superalgebra
still has a good description.

\subsection{}
The  basic classical Lie superalgebras were described by V.~Kac
in~\cite{kadv}. These Lie superalgebras are the 
closest to the ordinary simple Lie algebras: for instance, they can be 
described by their Cartan matrices. The even part of a basic classical Lie 
superalgebra $\fg$ is a reductive Lie algebra.  
The centre $\Zg$ is described
in~\cite{kch1},~\cite{s1},~\cite{bzv}; it is isomorphic to
a subalgebra of $\Sh^W$ and their fields
of fractions coincide (here $\fh$ is a Cartan subalgebra
of $\fg_0$ and $W$ is the Weyl group of $\fg_0$).

\subsection{}
Let $\fg=\fg_0\oplus\fg_1$ be a basic classical Lie superalgebra, $\Ug$ be its
universal enveloping superalgebra and $\Zg$ be the centre of $\Ug$.
Let $T$ be a special {\em ghost} element constructed in~\cite{ghost}---
see~\ref{TUg'}. We call a maximal ideal $\wdchi$ of $\Zg$
{\em strongly typical} if it does not contain $T^2$.

For a maximal ideal $\wdchi$ of $\Zg$ and
a $\fg$-module $\wdN$ set $\wdN_{\wdchi}:=\{v\in\wdN|\ \wdchi^rv=0,\
 \forall r>>0\}$. For a fixed $\wdchi\in\Max\Zg$,
denote by $\gr\wdC_r$ the category of graded $\fg$-modules $\wdN$ satisfying
$\wdchi^r\wdN=0$ and by $\gr\wdC_{\infty}$ the 
category of graded $\fg$-modules $\wdN$ satisfying
$\wdN_{\wdchi}=\wdN$.

Consider $\fg_0$ as a purely even Lie superalgebra. Denote
by $\cU$ its enveloping algebra and by $\cZ(\fg_0)$
the centre of $\cU$. For a maximal ideal $\chi$ of $\cZ(\fg_0)$
and a $\fg_0$-module $N$ define $N_{\chi}$ as
above. For a fixed $\chi\in\cZ(\fg_0)$, denote by $\gr\cC_r$ 
the category of graded $\fg_0$-modules $N$ satisfying
$\chi^r N=0$ and by $\gr\cC_{\infty}$ the category of
graded $\fg_0$-modules  $N$ satisfying $N_{\chi}=\wdN$. 
In Sect.~\ref{sectEquiv} we prove the following

\subsubsection{}
\begin{thm}{intr:EC}
For a strongly typical $\wdchi\in\Max\Zg$ the category
$\gr\wdC_{\infty}$  is equivalent
to the category $\gr\cC_{\infty}$
where $\chi$ is a suitable maximal ideal of $\cZ(\fg_0)$.
The equivalence is given by the functors
$$\wdN\mapsto \wdN_{\chi},\ \ \ \ 
N\mapsto (\Ind_{\fg_0}^{\fg} N)_{\wdchi}.$$
The restriction of these functors provide also
the equivalence of $\gr\wdC_r$ with $\gr\cC_r$.
\end{thm}

A maximal ideal $\chi\in\Max\cZ(\fg_0)$ is 
suitable in the sense of the above theorem iff
for some projective Verma $\fg$-module $\wdM$ satisfying
$\wdchi\wdM=0$, the $\fg_0$-module
$\wdM_{\chi}$ is Verma and, moreover,
one has $\wdM=\Ug \wdM_{\chi}$. Note  that if 
$\chi$ is a perfect mate for $\wdchi$
in the sense of~\cite{g}, then it fulfills
these conditions. In~\cite{g}, Sect.8 we prove the existence of a perfect
mate $\chi$ for each strongly typical $\wdchi\in\Max\Zg$.
This allows one to describe the category $\gr\wdC_{\infty}$
corresponding to any strongly typical $\wdchi$ in terms
of~\Thm{intr:EC}.

The proof is based on the theorem stating that for
a strongly typical $\wdchi$ the
annihilator of a Verma $\fg$-module $\wdM$
satisfying $\wdchi\wdM=0$ is equal
to $\Ug\wdchi$--- see~\cite{g}.

An analogue of~\Thm{intr:EC} for the basic classical
Lie superalgebras of type I is proven in~\cite{psI}; 
for ``generic'' $\wdchi$ in type II case, it is
proven in~\cite{p}. Both proofs are based on
a realization of the categories $\gr\wdC_1, \gr\cC_1$ as
categories of $D$-modules on the corresponding flag varieties.

\subsection{}
\Thm{intr:EC} suggests a study of categories of representations
of $\fg$ corresponding to a non strongly typical character $\wdchi$.
Call such $\wdchi$ {\em weakly atypical}.
If $\wdchi$ is not typical, the category $\wdC_{\infty}$
has a very complicated structure--- for instance,
it contains non-trivial extensions of
non-isomorphic finite dimensional modules.
An ``intermediate'' case, when $\wdchi$ is  typical
but not strongly typical (it is possible only
for the types $B(m,n), G(3)$), seems to
be less complicated. We consider the case when
$\fg={\frak {osp}}(1,2l)$ ($B(0,2l)$) and the maximal ideal 
is a ``generic'' weakly atypical 
(see the condition~(\ref{condbeta})).
This means that the highest weight $\lambda$ of
a highest weight module annihilated by $\fm$ belongs to exactly one of 
the hyperplanes $S_{\beta}, \beta\in\Delta^+_1$
where $S_{\beta}:=\{\mu\in\fh^*|\ (\mu+\rho,\beta)=0\}$.

It is convenient  to substitute
the maximal ideal $\fm$ of $\Zg$ by the maximal ideal $\wdchi$
of the ``ghost centre'' $\tZg$  (see~\cite{ghost}
for definition) which contains $\fm$.  For $\fg={\frak {osp}}(1,2l)$
one has $\tZg=\Zg\oplus T\Zg$ and the ideal $\wdchi:=\fm\oplus T\Zg$
is a maximal ideal of $\tZg$. Define the categories
$\gr\wdC_r$ and $\gr\wdC_{\infty}$ for such $\wdchi\in\Max \tZg$
in the same way as it was done above for $\wdchi\in\Max \Zg$.
For fixed $\chi$, denote by $\cC_r$ the category of non-graded 
$\fg_0$-modules $N$ satisfying
$\chi^r N=0$ and by $\cC_{\infty}$  the full
subcategory of the category of non-graded $\fg_0$-modules  
consisting of the modules $N$ satisfying $N_{\chi}=\wdN$.

Under the above assumption on $\fm$, we prove
in Sect.~\ref{osp12l} the following

\subsubsection{}
\begin{thm}{intr:osp2}
The category $\gr\wdC_{\infty}$  is equivalent
to the category $\cC_{\infty}$
where $\chi$ is a suitable maximal ideal of $\cZ(\fg_0)$.
The equivalence is given by the functors
$$\wdN\mapsto \wdN_{\chi}\cap\wdN_0,\ \ \ \ 
N\mapsto (\Ind_{\fg_0}^{\fg} N)_{\wdchi}$$
where the grading on $\Ind_{\fg_0}^{\fg} N$
is determined by the assuming $N$ to be even.

The restriction of these functors provide also
the equivalence of $\gr\wdC_r$ with $\cC_r$.
\end{thm}

A maximal ideal $\chi\in\Max\cZ(\fg_0)$ is suitable in the 
sense of the above theorem iff
for some projective Verma $\fg$-module $\wdM$ satisfying
$\wdchi\wdM=0$, the $\fg_0$-module
$M:=\wdM_{\chi}\cap\wdM_0$ is Verma and $\wdM=\Ug M$. 
In~\ref{cnstr} we construct such a suitable ideal $\chi$ 
for each $\fm\in\Max\Zg$ satisfying the above assumption.

The proof is based on the  theorem stating that
the annihilator of a Verma $\fg$-module $\wdM$ satisfying
$\fm\wdM=0$ is equal to $\Ug\wdchi$ (see~\cite{gl2}, 6.2)
and on~\Prop{propAB} which describes the locally finite part $F(\wdM,\wdM)$
of endomorphisms of a Verma module $\wdM$ through the image
of $\Ug$ in it.

\subsection{Acknowledgments} The results of this paper were obtained when 
the author was a visitor at 
MSRI and at Max-Planck Institut f\"ur Mathematik at Bonn.
I express my gratitude to these institutions for the hospitality and excellent 
working conditions. I wish to thank V.~Serganova and I.~Penkov for  
helpful discussions.

\section{Preliminaries}
Everywhere in the paper  ${\frak g}=\fg_0\oplus\fg_1$ 
denotes one (unless otherwise
specified, an arbitrary one) of
the basic classical complex Lie superalgebras
${\frak {gl}}(m,n),\ {\frak {sl}}(m,n),\
{\frak {osp}}(m,n),\ {\frak {psl}}(n,n)$.
Each of these Lie superalgebras possesses the following properties:
it admits a ${\frak g}$-invariant
bilinear form which is non-degenerate on $[\fg,\fg]$
and the even part
${\frak g}_0$ is a reductive Lie algebra.

\subsection{Conventions}
In this paper the ground field is ${\Bbb C}$.
We denote by ${\Bbb N}^+$ the set of positive integers.
If $A$ is an algebra, $N$ is an $A$-module and
$X,Y$ are subsets of $A$ and $N$ respectively, 
we denote by $XY$ the submodule spanned
by the products $xy$ where $x\in X,y\in Y$.

For a ${\Bbb Z}_2$-homogeneous
element $u$ of a superalgebra denote by $d(u)$ its 
${\Bbb Z}_2$-degree. In all formulae where this notation is used,
$u$ is assumed to be ${\Bbb Z}_2$-homogeneous.

For a Lie superalgebra $\fm$ denote by $\cU(\fm)$ its universal enveloping
algebra and by $\cS(\fm)$ its symmetric algebra.
All modules in the text are assumed to be left modules
unless otherwise specified.
An $\fm$-module $N$ is called locally
finite if $\dim \cU(\fm)v<\infty$ for all $v\in N$.
Set $\Ug:=\cU(\fg)$ and $\cU:=\cU(\fg_0)$.

The symbol $\wdV$ (resp., $V$) is always used for
a simple $\fg$ (resp., $\fg_0$) module
and the symbol $\wdM$ (resp., $M$) for
a Verma $\fg$ (resp., $\fg_0$) module.

\subsection{}
Fix a triangular decomposition $\fg=\fn^-\oplus\fh\oplus\fn^+$---
see~\cite{psg} for definition. Denote by $\Delta$ the set of
all non-zero roots of $\fg$ and by $\Delta_0^+$ (resp., $\Delta_1^+$)
the set of non-zero positive even (resp., odd) roots of $\fg$.
Set $\ol {\Delta}_1^+:=\{\beta\in\Delta^+_1|\ 2\beta\not\in\Delta_0^+\}$.
and
$$\rho:={1\over 2}(\sum_{\alpha\in\Delta^+_0}\alpha-
\sum_{\beta\in\Delta^+_1}\beta).$$

Let $W$ be the Weyl group of $\fg_0$. For $w\in W,\mu\in\fh^*$ set
$$w.\mu=w(\mu+\rho)-\rho.$$
Denote by $(-,-)$ a $\fg$-invariant bilinear form
on $\fg$ which is non-degenerate on $[\fg,\fg]$
and also the induced $W$-invariant bilinear form on $\fh^*$.
One has  $\ol{\Delta}_1^+=\{\beta\in\Delta^+_1|\ (\beta,\beta)=0\}$.

\subsection{}
\label{TUg'}
Define the adjoint action of $\fg$ on $\Ug$ by setting
$$(\ad g)u=gu-(-1)^{d(g)d(u)}ug, \ \forall g\in\fg,u\in\Ug$$
By default, the action of $\fg$ on $\Ug$ is assumed to be the adjoint action.
The centre $\Zg$ of $\Ug$ is equal to $\Ug^{\ad\fg}$.

Define the twisted adjoint action of $\fg$ on $\Ug$ by setting
$$(\ad'g)u=gu-(-1)^{d(g)(d(u)+1)}ug, \ \forall g\in\fg,u\in\Ug.$$
The anticentre $\Ag:=\Ug^{\ad'\fg}$ contains an element
$T$ defined in~\cite{ghost}. This is a unique element
of $\Ag$ satisfying
$\cP(T)(\lambda)=\prod_{\beta\in\Delta_1^+}(\beta,\lambda+\rho)$
for any $\lambda\in\fh^*$; here
$\cP:\Ug\to\Sh$ is the Harish-Chandra projection. 
The element $T$ is even.
Since $T$ belongs to $\Ag$, it commutes with the even
elements of $\Ug$  and anticommutes with the odd ones; in particular,
$T^2\in \Zg$. 

The Harish-Chandra projection $\cP$ provides 
a monomorphism $\iota:\Zg\to\Sh^{W.}$. The image of $\iota$ is described 
in~\cite{kch1},~\cite{s1},~\cite{bzv}; $\iota$
is bijective iff $\fg={\frak {osp}}(1,2l)$.
The centre $\Zg$ contains an element $Q$ such that
$$\cP(Q)(\lambda)= \prod_{\beta\in\ol{\Delta}_1^+}(\beta,\lambda+\rho).$$
The localized algebra $\Zg[T^{-2}]$
is isomorphic to a localization of a polynomial algebra
$\Sh^{W}$ by $\cP(Q)$---see~\cite{kch}. 
The localized algebra $\Zg[T^{-2}]$
is isomorphic to a localization of a polynomial algebra
$\Sh^{W}$.

\subsection{The category $\wdO$ and Verma modules}
\label{wdO}
Denote by $\cO$ the full subcategory of the category
of $\fg_0$-modules consisting of finitely generated
$\fh$-diagonalizable $\fg_0$-modules  which are $\fn^+_0$-locally
finite. Denote by $\wdO$ the similarly defined category
of $\fg$-modules. 
Since $\cU(\fn^+)$ is finite over $\cU(\fn_0^+)$, 
a $\fg$-module $N$ belongs $\wdO$ iff as a 
$\fg_0$-module $N$ belongs to $\cO$. In particular, 
any  module of category $\wdO$ has a finite length.

\subsubsection{}
\label{Vsmod}
For $\lambda\in {\frak h}^*$ denote by ${\Bbb C}_{\lambda}$
a one-dimensional ${\frak b}$-module such that
${\frak n}^+v=0$ and $hv=\lambda(h)v$ for any $h\in {\frak h},v\in  
{\Bbb C}_{\lambda}$. Define a Verma module $\wdM(\lambda)$ by setting
$$\wdM(\lambda):={\Ug}\otimes_{\cU({\frak b})}
{\Bbb C}_{\lambda}.$$
A Verma $\fg_0$-module $M(\lambda)$ is defined similarly.

\subsubsection{}
\label{W()}
A maximal ideal $\wdchi\in\Max\Zg$ is called {\em typical}
if it does not contain $Q$ defined in~\ref{TUg'}.
If $\wdchi\in\Max\Zg$ is typical then the set 
$$W(\wdchi):=\{\lambda\in\fh^*|\ \wdchi\wdM(\lambda)=0\}$$
forms a single $W.$-orbit. For $\fg={\frak {osp}}(1,2l)$
all $\wdchi\in\Max\Zg$ are typical.

Call a $\fg$-module $N$ typical if $\Ann_{\Zg} N$
is a typical maximal ideal of $\Zg$.
Define the partial order on $\fh^*$ by setting
$\mu\geq\nu\  \Longleftrightarrow\ 
(\mu-\nu)\in \sum_{\alpha\in \Delta^+}{\Bbb N}\alpha$. 
A typical Verma module $\wdM(\lambda)$ is projective (resp., simple)
in $\wdO$ if $\lambda$ is maximal (resp., minimal) in $W.\lambda$---see, 
for instance,~\cite{g}, 2.5.3.
In particular, for a typical $\wdchi\in\Max\Zg$
there exists a projective Verma module
$\wdM$ satisfying $\wdchi\wdM=0$.

A maximal ideal $\wdchi\in\Max\Zg$ is called {\em strongly typical}
if $T^2\not\in \wdchi$. A strongly typical central character
is typical.  Call a $\fg$-module $N$ strongly typical if $\Ann_{\Zg} N$
is a strongly typical maximal ideal of $\Zg$.

If $\fg$ is not of types $B(m,n), G(3)$
then $\Delta^+_1=\ol{\Delta}^+_1$ and so 
the notion of typical and strongly typical module coincide.

\subsection{}
\label{F}
Throughout the paper we shall write ``$\ad\fg$-module'' instead
``$\fg$-module with respect to the adjoint action''.
For any $\fg$-modules $N_1,N_2$ view 
$\Hom(N_1,N_2):=\Hom_{\Bbb C}(N_1,N_2)$ as a $\fg$-module
with respect to the adjoint action and
denote by $F(N_1,N_2)$ the locally finite part of 
the $\ad\fg$-module $\Hom(N_1,N_2)$.
Similarly for  any $\fg_0$-modules $N_1,N_2$
denote by $F(N_1,N_2)$ the locally finite part of 
the $\ad\fg_0$-module $\Hom(N_1,N_2)$.
Notice that for any $\fg$-module $N$ its locally finite part
coincides with its $\ad\fg_0$-locally 
finite part, since $\Ug$ is a finite extension
of $\cU$.

\subsubsection{}
\label{FNE}
Let $N$ be a $\fg_0$-module and $E$ be a finite dimensional
$\fg_0$-module. It is easy to check that 
$A:=F(N,N),B:=F(N\otimes E,N\otimes E)$ are subalgebras
of $\Hom(N,N)$ and  $\Hom (N\otimes E,N\otimes E)$ respectively.
The algebra $A$ acts on $F(N,N\otimes E)$ from the right
and the algebra $B$ acts on $F(N,N\otimes E)$ from the left;
these actions commute.
We claim that $F(N,N\otimes E)$ is a free right $A$-module
whose rank is equal to the dimension of $E$ and, moreover,
$\End_A(F(N,N\otimes E))=B$. 

Indeed, consider the  map
$\iota: F(N,N)\otimes E\to \Hom(N,N\otimes E)$ given
by $\iota(f\otimes v)(n)=f(n)\otimes v$ for any $f\in F(N,N), v\in E$
and the map $\iota': F(N,N\otimes E)\to \Hom(N,N)\otimes E$
given by $\psi\mapsto \sum_i p_i\circ\psi\otimes e_i$ where
$\{e_i\}$ is a basis of $E$ and $p_i: N\otimes E\to N$
are given by $p_i(n\otimes e_j):=\delta_{i,j} n\otimes e_j$.
One can easily sees that $\im\iota\subseteq F(N,N\otimes E),
\im\iota'\subseteq F(N,N)\otimes E$ and that
$\iota'\circ\iota=\id, \ \iota\circ\iota'=\id$.  
Thus $F(N,N\otimes E)\cong F(N,N)\otimes E$ is a free right $A$-module
whose rank is equal to the dimension of $E$.
Similarly, the map $\iota'':A\otimes \End_{\Bbb C} (E)\to B$
given by $\iota''(f\otimes \phi)(n\otimes v):=f(n)\otimes \phi(v)$
for any $f\in A=F(N,N), n\in N,v\in E,\phi\in \End_{\Bbb C} (E)$
is bijective. For any $f,f'\in A, v\in E,\phi\in \End_{\Bbb C} (E)$
one has $\iota''(f\otimes \phi)(\iota(f'\otimes v))=
\iota(ff'\otimes \phi(v))$. This implies $\End_A(F(N,N\otimes E))=B$
since the $A$-module $F(N,N\otimes E)$ is freely generated by
the elements $\iota(1\otimes e_i)$ (here $1$ is the unit of $A$).

\subsection{}
\label{ann}
Let $M$ be a Verma $\fg$-module.
By Duflo's theorem $\Ann M=\cU\Ann_{\cZ(\fg_0)} M$---see~\cite{d}.
By~\cite{j16}, 6.4 the natural map $\cU/(\Ann M)\to F(M,M)$
is bijective.

Let $\wdM$ be a strongly typical Verma $\fg$-module. Then
$\Ann\wdM=\Ug\Ann_{\Zg}\wdM$ 
and the natural map $\Ug/(\Ann\wdM)\to F(\wdM,\wdM)$
is bijective (see~\cite{g}, 9.4,9.5).

\subsection{}
For a $\fg$-module $\wdN$ and a maximal ideal $\wdchi\in\Zg$ 
set 
$$\wdN_{\wdchi}:=\{v\in \wdN|\ \wdchi^rv=0, \ r>>0\}.$$
We say that a $\Ug$-module $\wdN$ has {\em a finite support }
$\supp_{\Zg}\wdN=\{\wdchi_1,\ldots,\wdchi_k\}$ if
for any $v\in\wdN$ there exist $r_1,\ldots,r_k\in {\Bbb N}^+$
such that $\prod_i\wdchi_i^{r_i}v=0$. In this case,
$$\wdN=\oplus_i \wdN_{\wdchi_i}$$
and each $\wdN_{\wdchi_i}$ is canonically isomorphic to the localization
of the module $\wdN$ at $\wdchi_i$.
If $\wdN$ has a finite support and $0\to\wdN'\to\wdN\to \wdN''\to 0$
is an exact sequence then, for any $\wdchi'\in\Max\Zg$, the sequence
$0\to\wdN'_{\wdchi'}\to\wdN_{\wdchi'}\to \wdN''_{\wdchi'}\to 0$
is also exact.

We adopt the similar notation for $\cZ(\fg_0)$ and $\cU$-modules.
For a graded $\fg$-module $\wdN$ and a maximal ideal $\chi$ of
$\cZ(\fg_0)$ we set
$$\wdN_{\chi;i}:=\wdN_{\chi}\cap \wdN_i$$
for $i=0,1$.

\subsubsection{}
\label{finsupp}
Let $\wdchi$ be a strongly typical central character,
$\wdM$ be such that $\wdchi\wdM=0$. Then 
$$\Ug\wdchi\cap \cZ(\fg_0)=\Ann\wdM\cap \cZ(\fg_0)=
\prod_{\chi\in\supp_{\cZ(\fg_0)}\wdM} \chi^{r(\chi)}$$
where $r(\chi)$ is the minimal $r$ such that $\chi^r\wdM_{\chi}=0$.
In particular, if $\chi\in\Max\cZ(\fg_0)$
is such that $\wdM_{\chi}$ is a Verma $\fg_0$-module
then $r(\chi)=1$. Moreover, any $\fg$-module $\wdN$ satisfying
$\wdchi\wdN=0$ has a finite support in $\Zg$ (which is a subset of 
$\supp_{\cZ(\fg_0)}\wdM$) and one has $\chi\wdN_{\chi}=0$ if  
$\wdM_{\chi}$ is a Verma $\fg_0$-module.

\subsubsection{}
\label{prfmate}
For a strongly typical $\wdchi\in\Max\Zg$ call
$\chi\in\Max\cZ(\fg_0)$ {\em a perfect mate} if 

(i) For any Verma $\fg$-module $\wdM$
annihilated by $\wdchi$, the $\fg_0$-module
$\wdM_{\chi}$ is Verma.

(ii) For any non-trivial $\fg$-module $\wdN$ annihilated by $\wdchi$,
the $\fg_0$-module $\wdN_{\chi}$ is non-trivial.

In~\cite{g}, Sect.8 we describe a perfect mate for
each strongly typical $\wdchi\in\Max\Zg$.

\subsection{}
\label{indcoind}
For a graded $\fg_0$-module $L$ 
denote by $\Ind_{\fg_0}^{\fg}L$ the vector space
$\Ug\otimes_{\cU} L$ (here $\Ug$ is considered as a right
$\cU$-module and a left $\Ug$-module through the multiplication)
equipped with the natural structure of a left graded
$\Ug$-module.  
Denote by $\Coind_{\fg_0}^{\fg}L$ the vector space
$\Hom_{\cU}(\Ug,L)$ (here $\Ug$ is considered as a left
$\cU$-module) equipped with the following structure of a left graded
$\Ug$-module: $(uf)(u'):=f(u'u)$ for any 
$f\in \Hom_{\cU}(\Ug,L),\ 
u,u'\in \Ug$. For a graded $\fg$-module $\wdN$ and a graded
$\fg_0$-module $L$ one has the canonical bijections
\begin{equation}
\label{isoind}
\begin{array}{cl}
\Hom_{\fg_0}(\wdN,L)\iso
\Hom_{\fg}(\wdN,\Coind_{\fg_0}^{\fg}L),& \\
\Hom_{\fg_0}(L, \wdN)\iso
\Hom_{\fg}(\Ind_{\fg_0}^{\fg}L, \wdN). &
\end{array}\end{equation}
By~\cite{bf}, $\Ind_{\fg_0}^{\fg}L\cong\Coind_{\fg_0}^{\fg}L$
for any graded $\fg_0$-module $L$.

The same formulae define non-graded versions of $\Ind_{\fg_0}^{\fg}L$ and 
$\Coind_{\fg_0}^{\fg}L$. The same canonical bijections~(\ref{isoind})  
take place.

\section{Equivalence of Categories for a strongly typical central character}
\label{sectEquiv}
In this section we prove that the category
of $\Ug$-modules $\wdN$ satisfying $\wdN=\wdN_{\wdchi}$
is equivalent to the category
of $\cU$-modules $N$ satisfying $N=N_{\chi}$
provided $\wdchi\in \Max\Zg$ is strongly typical and 
$\chi\in\Max\cZ(\fg_0)$ is its perfect mate---see~\ref{prfmate}.
The corresponding quasi-inverse functors are 
given by $\wdN\mapsto \wdN_{\chi}$ and 
by $N\mapsto (\Ind_{\fg_0}^{\fg} N)_{\wdchi}$.
There are two versions of this result: graded and non-graded.
We consider first a non-graded version. A graded version is easily deduced
from the non-graded one.

\subsection{Notation}
Take a strongly typical $\wdchi\in\Max\Zg$ and let $\chi\in\Max\cZ(\fg_0)$.
Denote by $\wdC_r$ ($r\in {\Bbb N}^+$) the category
of non-graded $\Ug$-modules $\wdN$ satisfying $\wdchi^r\wdN=0$
and by $\wdC_{\infty}$ the category
of non-graded $\Ug$-modules $\wdN$ satisfying $\wdN=\wdN_{\wdchi}$.
Similarly, let $\cC_r$ be the category
of $\cU$-modules $N$ satisfying $\chi^r N=0$
and $\cC_{\infty}$ be the category
of $\cU$-modules $N$ satisfying $N=N_{\wdchi}$.
Evidently $\wdC_r$ is a full subcategory of $\wdC_{r+1}$
and any module in $\wdC_{\infty}$ 
is a direct limit of modules belonging to $\wdC_r$ for $r\to\infty$.

\subsubsection{}
\label{gen(b)}
Let $M$ be a Verma $\fg_0$-module such that $\chi M=0$.
We shall use the following equality
$$\Ann_{\Ug} (\Ind_{\fg_0}^{\fg} M)=\Ug\chi$$
which follows from the fact that $\Ug$
is free over $\cU$.

Let us show that 
\begin{equation}
\label{suppind}
\supp_{\Zg}\Ind_{\fg_0}^{\fg}N\subseteq\supp_{\Zg}\Ind_{\fg_0}^{\fg}M
\end{equation}
for any $N\in \cC_{\infty}$.
Indeed, for any $N\in\cC_1$ 
$$\Ann_{\Zg}(\Ind_{\fg_0}^{\fg}N)\supseteq
\Ug\chi\cap\Zg=\Ann_{\Zg}(\Ind_{\fg_0}^{\fg}M)$$
that implies the inclusion~(\ref{suppind}) (for $N\in\cC_1$).
Any $N\in\cC_r$ admits a finite filtration
with the factors belonging to $\cC_1$
and so  the inclusion~(\ref{suppind}) holds
for such $N$. To deduce~(\ref{suppind})
for any $N\in \cC_{\infty}$, observe that for
any $v\in \Ind_{\fg_0}^{\fg}N$
there exists a finitely generated submodule $N'$ of $N$
such that $\Ind_{\fg_0}^{\fg}N'$ contains $v$.
Since $N'$ is finitely generated, it lies in $\cC_r$ for
a suitable $r\in {\Bbb N}^+$. This implies  the inclusion~(\ref{suppind}).
 
Hence   for any $N\in \cC_{\infty}$
\begin{equation}\label{ind=oplus}
\Ind_{\fg_0}^{\fg}N=\oplus_{\wdchi'\in \supp _{\Zg}\Ind_{\fg_0}^{\fg}M}
(\Ind_{\fg_0}^{\fg}N)_{\wdchi'}
\end{equation}

Similarly,  choose $\wdM$ such that $\Ug\wdchi=\Ann \wdM$. Then
for any $\wdN\in\wdC_{\infty}$ 
\begin{equation}\label{wdn=oplus}
\wdN=\oplus_{\chi'\in \supp_{\cZ(\fg_0)}\wdM}
\wdN_{\chi'}.
\end{equation}

\subsubsection{}
Consider the map $f:\Ug/(\Ug\chi)\to F(M,\Ind_{\fg_0}^{\fg}M)$
induced by the natural map $\Ug\to F(M,\Ind_{\fg_0}^{\fg}M)$
given by $u\mapsto(m\mapsto u\otimes m)$.
We claim $f$ is an isomorphism
of left $\Ug$ and right $\cU$-modules. Since
$\Ug$ is free over $\cU$ and $\Ann_{\cU} M=\cU\chi$, the map
$f$ is injective.
Kostant's Separation Theorem (see~\cite{ko}) 
states the existence of an $\ad\fg_0$-submodule
$H'$ of $\cS(\fg_0)$ such that the multiplication map
provides an isomorphism $H'\otimes\cS(\fg_0)^{\fg_0}\iso\cS(\fg_0)$.
Then
$$\Sg=\Lambda\fg_1\otimes\cS(\fg_0)=(\Lambda\fg_1\otimes H')
\otimes\cS(\fg_0)^{\fg_0}.$$
Using~\cite{bl}, 5.4, it is easy to deduce the existence 
of an $\ad\fg_0$-submodule $H$ of $\Ug$ such that the multiplication
map provides the isomorphism $H\otimes\cZ(\fg_0)\to \Ug$.
Moreover, $H\cong \Lambda\fg_1\otimes H'$ as an $\ad\fg_0$-module.
By~\ref{ann}, $H'\cong  F(M,M)$ as $\ad\fg_0$-modules.
The following chain of $\ad\fg_0$-isomorphisms
$$\Ug/(\Ug\chi)\cong\Lambda\fg_1\otimes H'\cong
\Lambda\fg_1\otimes F(M,M)\cong  F(M,\Ind_{\fg_0}^{\fg}M)$$
implies the surjectivity of $f$. 

The bijectivity of $f$ gives the following useful formula
\begin{equation}
\label{formulaInd}
\Ind_{\fg_0}^{\fg}N= F(M,\Ind_{\fg_0}^{\fg}M)\otimes_{\cU}N.
\end{equation}
for any $N\in \cC_1$.

\subsection{}
Fix any $\chi\in\Max\cZ(\fg_0),\wdchi\in\Max\Zg$. Define
the functors $\Psi:\wdC_{\infty}\to\cC_{\infty}$ and
$\Phi: \cC_{\infty}\to\wdC_{\infty}$ by the formulae
$$
\Psi(\wdN)=\wdN_{\chi},\quad
\Phi(N)=(\Ind_{\fg_0}^{\fg} N)_{\wdchi}.$$

\subsubsection{}
\begin{lem}{adjfunct}
The functor $\Psi$ is left and right adjoint to $\Phi$.
\end{lem}
\begin{pf} Using~(\ref{ind=oplus}), one obtains 
for any $N\in\cC_{\infty},\wdN\in\wdC_{\infty}$ 
$$\begin{array}{rl}
\Hom_{\fg}(\Phi(N),\wdN)=\Hom_{\fg}((\Ind_{\fg_0}^{\fg} N)_{\wdchi},\wdN)&=
\Hom_{\fg}(\Ind_{\fg_0}^{\fg} N,\wdN)=\Hom_{\fg_0}(N,\wdN)\\
&=\Hom_{\fg_0}(N,\wdN_{\chi})=\Hom_{\fg_0}(N,\Psi(\wdN))
\end{array}$$
and also, using~(\ref{wdn=oplus}),
$$\begin{array}{rl}
\Hom_{\fg}(\wdN,\Phi(N))&=\Hom_{\fg}(\wdN,(\Ind_{\fg_0}^{\fg} N)_{\wdchi})=
\Hom_{\fg}(\wdN,\Ind_{\fg_0}^{\fg} N)\iso
\Hom_{\fg}(\wdN,\Coind_{\fg_0}^{\fg} N)\\
 &=
\Hom_{\fg_0}(\wdN,N)=\Hom_{\fg_0}(\wdN_{\chi},N)=\Hom_{\fg_0}(\Psi(\wdN),N)
\end{array}$$
where $=$ stands for the natural isomorphisms
and $\iso$ is induced by an isomorphism 
$\Ind_{\fg_0}^{\fg} N\cong \Coind_{\fg_0}^{\fg} N$ (see~\ref{indcoind}).
\end{pf}

\subsection{Conventions}
\label{conv}
Fix a strongly typical $\wdchi\in \Max\Zg$
and a Verma $\fg$-module $\wdM$ which is projective
in $\wdO$ and such that $\wdchi\wdM=0$.
Take $\chi\in\Max\cZ(\fg_0)$ such that
\begin{equation}
\label{propertyab}
\begin{array}{ll}
(a) & M:=\wdM_{\chi} \text{ is a Verma }\fg_0\text{-module}\\
(b) & \wdM=\Ug M.
\end{array}\end{equation}
For instance, one can choose $\chi$ to be a perfect mate
for $\wdchi$ (see~\ref{prfmate}).

Till the end of this section 
$\wdchi$, $\wdM$ and $\chi$ chosen as above
are assumed to be fixed.

\subsubsection{}
\begin{thm}{thmEquiv}
The functors 
$$\begin{array}{ll}
\Psi:\wdC_{\infty}\to\cC_{\infty}\ & 
\wdN\mapsto \wdN_{\chi},\\
\Phi: \cC_{\infty}\to\wdC_{\infty} \ &
N\mapsto (\Ind_{\fg_0}^{\fg} N)_{\wdchi}
\end{array}$$ 
are mutually quasi-inverse. Moreover, their restrictions
provide the equivalence of the categories
$\wdC_r$ and $\cC_r$ for any $r\in {\Bbb N}^+$.
\end{thm}

{\em Outline of the proof.}
We know by~\ref{adjfunct}  that the functors $\Psi,\Phi$ are adjoint.
In~\ref{redcase1} we reduce the required
assertion to the ``case $r=1$'' that is to the statement
that the restriction $\Psi_1,\Phi_1$ of the functors $\Psi,\Phi$
to the categories $\wdC_1,\cC_1$ 
provide an equivalence of the categories.
Not that~\ref{redcase1} does not use the condition~(\ref{propertyab}).

Observe that 
the inclusion $\Psi_1(\wdC_1)\subseteq\cC_1$ immediately
follows from~\ref{finsupp}, but it is not clear
apriori that $\Phi_1(\cC_1)\subseteq\wdC_1$.
To prove that $\Psi_1,\Phi_1$ provide an equivalence of the categories
we show, using~\Prop{equiv}, that the functor 
$\Phi': N\mapsto F(M,\wdM)\otimes_{\cU} N$ provides
an equivalence of the categories $\cC_1\to\wdC_1$.
It is easy to show that $\Psi_1$ is the left quasi-inverse
to $\Phi'$ and so it provides
an equivalence of the categories $\wdC_1\to\cC_1$.
Using this fact, we show that $\Phi(M)\cong\wdM$
and deduce that $\Phi_1$ is isomorphic to $\Phi'$.
This will complete the proof.

\subsection{Proof of Theorem~\ref{thmEquiv}}
\subsubsection{Reduction to the case $r=1$.}
\label{redcase1}
The formula~(\ref{ind=oplus}) implies
that $\Phi(N)$ is isomorphic to the localization
of $\Ind_{\fg_0}^{\fg}N$, considered as a $\Zg$-module,
at the maximal ideal $\wdchi$ 
(that is by the set $\Zg\setminus\wdchi$).

Similarly, the formula~(\ref{wdn=oplus}) implies that
$\Psi(\wdN)$  is isomorphic to the localization
of $\wdN$, considered as a $\cZ(\fg_0)$-module,
at the maximal ideal $\chi$.

Taking into account that the induction and 
the localization functors are exact and
commute with direct limits, one concludes that
the assertion of~\Thm{thmEquiv} is equivalent
to the statement that 
the restrictions $\Psi_1,\Phi_1$ of the functors $\Psi,\Phi$
to the categories $\wdC_1,\cC_1$ 
provide an equivalence of the categories.

\subsubsection{Case $r=1$}
\label{caser=1}
It remains to show that the restriction  $\Psi_1$ of $\Psi$ to
the subcategory $\wdC_1$ and
the restriction  $\Phi_1$ of $\Phi$ to
the subcategory $\cC_1$ provide an equivalence
of the categories $\wdC_1$ and $\cC_1$. 

We start with the following technical proposition.

\subsubsection{}
\begin{lem}{equiv}
Assume that $\wdN$ is a $\fg$-module, $N$ is
a $\fg_0$-direct summand of $\wdN$ and
$\wdN$ is a $\fg$-direct summand of $\Ind_{\fg_0}^{\fg}N$.
Then the $F(\wdN,\wdN)$-$F(N,N)$
bimodule $F(N,\wdN)$ provides a Morita equivalence between
the algebras $F(N,N)$ and $F(\wdN,\wdN)$. 
\end{lem}
\begin{pf}
Denote the algebra $F(N,N)$ by $A$.
For any $\fg$-module $X$ endow the vector space
$F(N,X)$ with the natural right $A$-module structure.
Since $N$ is a $\fg_0$-direct summand of $\wdN$,
$A$ as a right module over itself is a direct summand 
of $F(N,\wdN)$. To obtain the statement
one has to check that $F(N,\wdN)$ is a finitely generated
projective right $A$-module and
that $\End_A (F(N,\wdN))=F(\wdN,\wdN)$.

Denote $\Ind_{\fg_0}^{\fg}N$ by $I$. Recall that as a $\fg_0$-module
$I=\Lambda\fg_1\otimes N$. By~\ref{FNE}, $F(N,I)$
is a free $A$-module whose rank is equal to the dimension
of $\Lambda\fg_1$. Since $\wdN$ is a direct summand of $I$,
the $A$-module $F(N,\wdN)$ is a direct summand of 
$F(N,I)$. Hence $F(N,\wdN)$ is a finitely generated
projective right $A$-module. 
One has $\End_A (F(N,\wdN))=p\End_A (F(N,I))p$
where $p\in \End_A (F(N,I))$ is the idempotent with the image
$F(N,\wdN)$ corresponding to the decomposition
$F(N,I)=F(N,\wdN)\oplus F(N,G)$. By~\ref{FNE}, 
the left action of $F(I,I)$  on $F(N,I)$ induces an isomorphism
$F(I,I)\iso\End_A (F(N,I))$. Thus the left action of
$F(\wdN,\wdN)$ on $F(N,\wdN)$ induces an isomorphism
$F(\wdN,\wdN)\iso \End_A (F(N,\wdN))$.
The assertion follows.
\end{pf}

\subsubsection{}
Retain the notation of~\ref{conv}.
Consider a canonical map $\Ind_{\fg_0}^{\fg} M\to\wdM$. 
By the property (b) of~(\ref{propertyab}), this map
is surjective. Recall $\wdM$ is projective in $\wdO$.
Consequently, $\wdM$ is a $\fg$-direct summand of $\Ind_{\fg_0}^{\fg} M$
and so the pair $(\wdM,M)$ satisfies the assumptions of~\Lem{equiv}.
From~\ref{ann} 
it follows that $F(\wdM,\wdM)\cong \Ug_{\wdchi}$ and $F(M,M)\cong \cU_{\chi}$
as algebras. Taking into account~\Lem{equiv} one concludes
that the functor $\Phi':\cC_1\to \wdC_1$ defined
by 
$$\Phi'(N):= F(M,\wdM)\otimes_{\cU_{\chi}} N= F(M,\wdM)\otimes_{\cU} N$$ 
provides an equivalence of the categories.

\subsubsection{}
\label{laststep}
One has
$$\Psi_1\circ\Phi'(N)=F(M,\wdM_{\chi})\otimes_{\cU} N=
F(M,M)\otimes_{\cU} N=N$$
since $F(M,M)\cong\cU_{\chi}$.
Hence $\Psi_1$ is quasi-inverse to $\Phi'$ and it
provides an equivalence of the categories
$\wdC_1\to\cC_1$.

Let us verify  that the functors $\Phi_1$ and $\Phi'$
are isomorphic. Formula~(\ref{formulaInd}) implies that
\begin{equation}
\label{formulaPhi}
\Phi_1(N)=F(M,(\Ind_{\fg_0}^{\fg}M)_{\wdchi})\otimes_{\cU}N=
F(M,\Phi_1(M))\otimes_{\cU}N
\end{equation}
for any $N\in\cC_1$.
Thus to show that $\Phi_1(N)\cong\Phi'(N)$ for all $N\in\cC_1$
it is enough to verify only that 
\begin{equation}
\label{IndM}
(\Ind_{\fg_0}^{\fg}M)_{\wdchi}\cong\wdM.
\end{equation}
As we have shown above,  
$\wdM$ is a $\fg$-direct summand of $\Ind_{\fg_0}^{\fg} M$.
Thus it is enough to check that 
$$\Hom_{\fg}(\wdN,\Ind_{\fg_0}^{\fg}M)=0$$
for all simple $\wdN\in \wdC_1$ such that $\wdN\not\cong\Soc\wdM$
and that
$$\dim\Hom_{\fg}(\Soc \wdM,\Ind_{\fg_0}^{\fg}M)=1.$$
Recall that
$\Ind_{\fg_0}^{\fg}M\cong\Coind_{\fg_0}^{\fg}M$ (see~\ref{indcoind}) and so
$$\Hom_{\fg}(\wdN,\Ind_{\fg_0}^{\fg}M)=
\Hom_{\fg_0}(\wdN,M)=\Hom_{\fg_0}(\wdN_{\chi},M).$$
If $\wdN\in \wdC_1$ is simple, $\wdN_{\chi}=\Psi_1(\wdN)\in\cC_1$ is
also simple.
Taking into account that the $\fg_0$-socle of $M$ is
a simple Verma $\fg_0$-module, one concludes that
$\Hom_{\fg_0}(\wdN_{\chi},M)=0$ if $\wdN_{\chi}\not\cong\Soc M$
and $\dim\Hom_{\fg_0}(\wdN_{\chi},M)=1$ otherwise.
If $\Psi_1(\wdN)=\wdN_{\chi}\cong\Soc M$ then $\wdN\cong\Soc \wdM$
since $\Psi_1:\wdC_1\to\cC_1$ provides an equivalence of the categories
and $\Psi_1(\wdM)=M$. This proves~(\ref{IndM}). 

Hence the functors $\Phi_1$ and $\Phi'$
are isomorphic. This completes the proof of~\Thm{thmEquiv}.

\subsection{Remarks}
\subsubsection{}
\label{etc}
Recall that any $\wdN\in\wdC_{\infty}$ is generated by
$\wdN_{\chi}$ and so is a homomorphic image
of $\Ind_{\fg_0}^{\fg} \wdN_{\chi}$.
Consequently, the restrictions
of  the functors $\Psi,\Phi$ provide
an equivalence between $\wdC_{\infty}\cap\wdO$ 
and $\cC_{\infty}\cap\cO$. The same is true for
the categories of weight modules.
Also, one can easily see that the
Verma modules of $\wdC_{\infty}\cap\wdO$ 
correspond to the Verma modules of $\cC_{\infty}\cap\cO$.

\subsubsection{}
\label{anway}
Formula~(\ref{IndM}) which says that
$\Phi(M)\cong\wdM$ is crucial for the whole proof. 
If one had an independent proof of this formula,
one could deduce from it a straightforward proof of~\Thm{thmEquiv}.

\subsection{Graded case}
\label{grdcase}
Denote by $\gr\wdC_{\infty}$ (resp., $\gr\wdC_r$) the full category
of ${\Bbb Z}_2$-graded $\Ug$-modules $\wdN$
satisfying $\wdN_{\wdchi}=\wdN$ (resp., $\wdchi^r\wdN=0$). 
Consider $\cU$ 
as a purely even superalgebra and define similarly 
$\gr\cC_{\infty},\gr\cC_r$.
Denote by $\#$ the forgetful functors
$\gr\wdC_{\infty}\to \wdC_{\infty}$
and $\gr\cC_{\infty}\to \cC_{\infty}$.

\subsubsection{}
Evidently any $N'\in\cC_{\infty}$ is isomorphic
to $N^{\#}$ for some $N\in\gr\cC_{\infty}$---
for instance, one can consider $N'$ as a purely even (odd) module.
It turns out that the similar assertion holds for 
$\wdC_{\infty}$. Indeed,
since $\wdchi$ is strongly typical, it contains an
element $(T^2-c^2)$ for some non-zero $c\in{\Bbb C}$.
Denote by $\Ug_0$ the even part of $\Ug$; recall
that $T$ lies in the centre of $\Ug_0$.
Any $\wdN\in\wdC_{\infty}$ is a direct sum
of the $\Ug_0$-modules $N_{+}$, $N_-$ where
$$N_+:=\{v\in\wdM|\ (T-c)^rv=0,\ \forall r>>0\},\ \ \ 
N_-:=\{v\in\wdM|\ (T+c)^rv=0,\ \forall r>>0\}.$$
One can define a ${\Bbb Z}_2$-grading on $\wdN$
by putting $\wdN_0:=N_+$, $\wdN_1:=N_-$.
In such a way, one obtains a functor
$\wdC_{\infty}\to \gr\wdC_{\infty}$
which is left quasi-inverse to the functor $\#$.
This implies, in particular, that for any
irreducible graded $\Ug_{\wdchi}$-module $\wdN$, its image
$\wdN^{\#}$ remains irreducible (non-graded) module.

\subsubsection{}
\begin{thm}{thmgrEquiv}
The functors 
$$\begin{array}{ll}
\Psi^{\gr}:\gr\wdC_{\infty}\to\gr\cC_{\infty},
\ & \wdN\mapsto \wdN_{\chi},\\
\Phi^{\gr}:\gr\cC_{\infty}\to\gr\wdC_{\infty},
\ & N\mapsto (\Ind_{\fg_0}^{\fg} N)_{\wdchi}
\end{array}$$ 
are mutually quasi-inverse. Their restriction
provides an equivalence of the categories
$\gr\wdC_r$ and $\gr\cC_r$ for any $r\in {\Bbb N}^+$.
\end{thm}
\begin{pf}
Repeating the arguments  of~\Lem{adjfunct}, one shows that the functor
$\Psi^{\gr}$ is left and right adjoint to the functor $\Phi^{\gr}$.

The functors $\Psi,\Phi$ are quasi-inverse and so the canonical
homomorphism $\alpha_N: N\to\Psi\circ\Phi(N)$ is
an isomorphism for any $N\in\cC_{\infty}$.
For any $N\in\gr\cC_{\infty}$ denote by $\alpha_N$ the canonical homomorphism
$N\to\Psi^{\gr}\circ\Phi^{\gr}(N)$. 
One has $\#\circ\Psi^{\gr}=\Psi\circ\#$ and $\#\circ\Phi^{\gr}=\Phi\circ\#$.
Therefore $\#(\alpha_N)=\alpha_{N^{\#}}$.
Since $\alpha_{N^{\#}}$ is an isomorphism, $\alpha_N$  is also an isomorphism.
Similarly, the canonical
homomorphism $\beta_{\wdN}: \Phi^{\gr}\circ\Psi^{\gr}(\wdN)\to\wdN$ is
an isomorphism for any $\wdN\in\gr\wdC_{\infty}$.
Hence $\gr\Psi$ and $\gr\Phi$ are quasi-inverse.

The equality $\Psi^{\gr}(\gr\wdC_r)=\gr\cC_r$ follows from
the equality $\Psi(\gr\wdC_r)=\cC_r$. This implies the second claim.
\end{pf}

\section{An Example}
\label{osp12l}
Say that a $\Ug$-central character  is weakly atypical if
it is  not strongly typical. In this section we consider 
a ``generic''  weakly atypical central characters 
for $\fg:={\frak {osp}}(1,2l)$.
We show that, for such a character  $\fm\in\Max\Zg$,
the category $\gr\wdC_{\infty}$ of
graded $\fg$-modules $\wdN$ satisfying $\wdN_{\fm}=\wdN$
is equivalent to the category $\cC_{\infty}$
of $\fg_0$-modules $N$ satisfying $N_{\chi}=N$
for an appropriate $\chi$. We see that in this case
the category $\gr\wdC_{\infty}$ is ``twice smaller''
than one could expect in the strongly typical case.

Throughout this section all $\fg$-modules are assumed to be graded.

\subsection{}
\label{annosp12l}
The superalgebra ${\frak {osp}}(1,2l)$ does not have
isotropic roots: $(\beta,\beta)\not=0$ for any $\beta\in\Delta_1$.
As a result, ${\frak {osp}}(1,2l)$ has many features of
the simple Lie algebras. For instance, $\Ug$ is a domain,
its centre $\Zg$ is a polynomial algebra
and all finite dimensional representations are completely reducible.
In~\cite{mu}, Musson proved the existence  of a harmonic space $H$ in $\Ug$:
this is an $\ad\fg$-submodule of $\Ug$ such that the multiplication
map provides the isomorphism $H\otimes\Zg\iso\Ug$.

If $\wdM$ is strongly typical its annihilator is
a centrally generated ideal. 

Suppose that $\wdM$ is not strongly typical; this means that
$\fm:= \Ann_{\Zg}\wdM$ contains $T^2$. Then, by~\cite{gl2}, 6.2,
$$\Ann\wdM=\Ug(T\Zg+\fm)$$
and $\wdchi:=(T\Zg+\fm)$ is the maximal
ideal of the ``ghost centre'' $\tZg=\Zg+T\Zg$.
The ideal $\Ug\wdchi$ is primitive since a Verma module
$\wdM(\mu)$ is simple if $\mu\in W(\fm)$ is minimal.

Moreover, by~\cite{gl2}, 5.3, for any simple $\fg_0$-module $V$ 
one has
\begin{equation}
\label{12}
\dim \Hom_{\fg_0}(V,\Ug/(\Ann \wdM))={1\over 2} \dim\Hom_{\fg_0}(V, H)
\end{equation}
where $H$ is the harmonic space mentioned above.

\subsection{}
Retain the notation of~\ref{W()}. Recall that $\fm\in\Max\Zg$ 
is strongly typical iff the elements $ \lambda\in W(\fm)$ do not
belong to the hyperplanes 
$$S_{\beta}:=\{\mu\in\fh^*|\ (\mu+\rho,\beta)=0\}$$
for $\beta\in\Delta^+_1$. In this section we consider
$\fm\in\Max\Zg$  such that any $\lambda\in W(\fm)$ belongs
to exactly one hyperplane $S_{\beta}$ that is
\begin{equation}
\label{condbeta}
\begin{array}{lll}
 & \exists! \beta\in\Delta^+_1: & (\lambda+\rho,\beta)=0.\\ 
\end{array}
\end{equation}
In particular, $\fm$ is not strongly typical
and so $T^2\in\fm$.
Set
$$\wdchi:=T\Zg+\fm.$$

For any $\fg$-module $\wdN$
the vector space 
$$\wdN_{\wdchi}:=\{v\in\wdN|\ \wdchi^rv=0,\ \ r>>0\}$$
is a $\fg$-submodule of $\wdN$.
Let $\gr\wdC_r$  (resp., $\gr\wdC_{\infty}$) 
be the category
of graded $\fg$-modules $\wdN$ satisfying $\wdchi^r\wdN=0$
(resp., $\wdN_{\wdchi}=\wdN$). Observe
that $\wdC_1\subseteq \wdC'_1\subseteq \wdC_2$ where $\wdC'_1$
is the category 
of graded $\fg$-modules $\wdN$ satisfying $\fm\wdN=0$.

Fix a projective Verma module $\wdM$ annihilated
by $\wdchi$. Fix $\chi\in\Max\cZ(\fg_0)$ such that
\begin{equation}
\label{condab}
\begin{array}{ll}
(a) & M:=\wdM_{\chi;0} \text{ is a Verma }\fg_0\text{-module}\\
(b) & \wdM=\Ug M.
\end{array}\end{equation}
(the existence of $\chi$ satisfying~(\ref{condab}) will be shown 
in~\ref{semiperfect}).
Let $\cC_{\infty}$ (resp., $\cC_{r}$) be the category
of $\fg_0$-modules $N$ satisfying $N_{\chi}=N$
(resp., $\chi^rN=0$).

\subsection{}
\begin{thm}{thmosp12l}
The functors
$$\begin{array}{ll}
\Psi: \gr\wdC_{\infty}\to \cC_{\infty} & \ \ \wdN\mapsto (\wdN_0)_{\chi},\\
\Phi: \cC_{\infty}\to \gr\wdC_{\infty}& \ \ N\mapsto
 (\Ind_{\fg_0}^{\fg} N)_{\wdchi}\end{array}$$
where the grading on $\Ind_{\fg_0}^{\fg} N=\Ug\otimes_{\cU}N$
is determined by assuming $N$ to be even,
are quasi-inverse. Moreover their restrictions
to $\gr\wdC_r$ and $\cC_r$ respectively
provide an equivalence of these categories.
\end{thm}
\begin{pf}
First steps of the proof are the same as those in
the proof of~\Thm{thmEquiv}. 
Considering  $\Ind_{\fg_0}^{\fg} N$ as the graded module
with respect to the grading defined above one obtains
$$\Hom_{\fg}(\Ind_{\fg_0}^{\fg} N,\wdN)=\Hom_{\fg_0}(N,\wdN_0)$$
for any graded $\fg$-module $\wdN$. Repeating~\ref{adjfunct}, one
concludes that the functors $\Psi,\Phi$ are adjoint.
Repeating~\ref{redcase1}, one reduces the required
assertion to ``the case $r=1$'' that is to the statement
that the restrictions $\Psi_1,\Phi_1$ of the functors $\Psi,\Phi$
to the categories $\gr\wdC_1,\cC_1$ 
provide an equivalence of the categories.

Evidently, $M$ is a $\fg_0$-direct summand of $\wdM$.
The embedding $M\to\wdM$ gives rise
to a non-zero $\fg$-map $\Ind_{\fg_0}^{\fg}M\to\wdM$.
Using $\wdM=\Ug M$ and the projectivity of $\wdM$, one concludes that $\wdM$
is a $\fg$-direct summand of $\Ind_{\fg_0}^{\fg}M$.
Denote by $B$ the algebra $F(\wdM,\wdM)$ and
by $B$-$\Mod$ the category of $B$-modules.
Taking into account~\Lem{equiv}, one concludes that the functor
$\Phi'':\cC_1\to B$-$\Mod$ defined
by 
$$\Phi''(N):= F(M,\wdM)\otimes_{\cU} N$$ 
provides an equivalence of the categories.

The category $\gr\wdC_1$ is the category of 
graded modules over the superalgebra $A:=\Ug/(\Ug\wdchi)$.
In~\Prop{propAB} below we prove that $B=A\oplus A\theta$ where
$\theta\in F(\wdM,\wdM)$ is given by
$$\theta(v)=(-1)^{d(v)}v.$$ 
In particular, $\theta^2=1$ and 
$\theta a=(-1)^{d(a)} a\theta$ for any $a\in A$.
Define the functor $\Gr: B$-$\Mod\to\gr\wdC_1$ as follows:
for any $X\in B$-$\Mod$ the $A$-module structure on $\Gr(X)$ is given
by ``the restriction of scalars'' and the
grading is given by
$$\Gr(X)_0:=\{v\in X|\ \theta v=v\},\ \ 
\Gr(X)_1:=\{v\in X|\ \theta v=-v\}.$$
It is easy to see that $\Gr$ provides an equivalence
of the categories  $B$-$\Mod$ and $\gr\wdC_1$.
Hence $\Phi':=\Gr\circ\Phi''$ provides an equivalence
of the categories $\cC_1$ and $\gr\wdC_1$.

The functor $\Phi': \cC_1\to\gr\wdC_1$
is given by the formula $N\mapsto F(M,\wdM)\otimes_{\cU} N$
where the grading on  $F(M,\wdM)$ is defined by
$F(M,\wdM)_i:=F(M,\wdM_i)$ ($i=0,1$).

The last step of the proof repeats~\ref{laststep}.
Indeed, 
$$\Psi_1\circ\Phi'(N)=F(M,(\wdM_{\chi;0})\otimes_{\cU} N=
F(M,M)\otimes_{\cU} N=N$$
since $F(M,M)\cong\cU_{\chi}$.
Hence $\Psi_1$ is quasi-inverse to $\Phi'$ and it
provides an equivalence of the categories
$\wdC_1\to\cC_1$. Now, in order 
to verify that the functors $\Phi_1$ and $\Phi'$
are isomorphic, one can simply repeat~\ref{laststep}
for the graded modules.
\end{pf}

\subsection{}
\label{etc1}
The restrictions of the functors $\Psi,\Phi$ provide
an equivalence between $\wdC_{\infty}\cap\wdO$ 
and $\cC_{\infty}\cap\cO$ and between
the corresponding categories of weight modules.
The Verma modules of $\wdC_{\infty}\cap\wdO$ 
correspond to the Verma modules of $\cC_{\infty}\cap\cO$.

\subsubsection{}
The category $\gr\wdC_{\infty}$ has a canonical involution
$\Pi$ given by
$$(\Pi\wdN)_0=\wdN_1,\ \ (\Pi\wdN)_1=\wdN_0.$$
This leads to an interesting involution $\Pi'=\Psi\circ\Pi\circ \Phi$ on the
category $\cC_{\infty}$. One has
$$\begin{array}{rl}
N&=((\Ind_{\fg_0}^{\fg} N)_{\wdchi})_{\chi;0},\\
\Pi'(N)&=((\Ind_{\fg_0}^{\fg} N)_{\wdchi})_{\chi;1}
\end{array}$$

\Thm{thmosp12l} implies that 
$M':=\wdM_{\chi;1}=\Psi(\Pi(\wdM))$ is a Verma $\fg_0$-module.
For  $N\in\cC_1$ one has
$\Phi(N)\cong\Phi'(N)$ and so
$$\Pi'(N)=F(M,M')\otimes_{\cU} N.$$

\subsection{}
\label{semiperfect}
It remains to prove~\Prop{propAB} and to find  
$\chi\in\Max\cZ(\fg_0)$ satisfying~(\ref{condab}).
To achieve these goals we proceed in a manner similar 
to~\cite{g}, Sect.8.

Call $\chi\in\Max\cZ(\fg_0)$ {\em a mate} for $\wdchi$
if for a Verma module $\wdM$ annihilated
by $\wdchi$ the $\fg_0$-modules
$\wdM_{\chi;0},\wdM_{\chi;1}$ are Verma modules.
Call $\chi\in\Max\cZ(\fg_0)$ {\em a perfect mate} for $\wdchi$
if it is a mate and for a simple highest weight module 
$\wdV(\lambda)$ annihilated by $\wdchi$ the $\fg_0$-module
$\wdV(\lambda)_{\chi;0}$ is non-zero.

In this subsection we construct a perfect mate $\chi$ for $\wdchi$ 
satisfying~(\ref{condbeta}).

\subsubsection{}
\label{Gamma}
The root system of $\fg$ takes form
$$\Delta_1^+=\{\sigma_i\}_{1\leq i\leq l},\ \ \ \ 
\Delta_0^+=\{\sigma_i\pm\sigma_j;2\sigma_i\}_{1\leq i<j\leq l}$$
and $(\sigma_i,\sigma_j)=\delta_{i,j}$. The Weyl group
$W$ acts  on $\{\sigma_i\}_1^l$ by the signed permutations.

Set 
$$\begin{array}{rl}
\Gamma :&=\left\{\displaystyle\sum_1^l r_i\sigma_i|\ r_i\in\{0,1\}\right\},\\
\Gamma_0:&=\left\{\displaystyle\sum_1^l r_i\sigma_i|\ r_i\in\{0,1\}, 
\displaystyle\sum_1^l r_i\text{ is even}\right\},\\
\Gamma_1:&=\left\{\displaystyle\sum_1^l r_i\sigma_i|\ r_i\in\{0,1\}, 
\displaystyle\sum_1^l r_i\text{ is odd}\right\}.
\end{array}$$
Define the action of the Weyl group $W$ on $\Gamma$
by setting
$$w_*\gamma=w(\gamma-\rho_1)+\rho_1.$$

Take an arbitrary $\lambda\in \fh^*$ and 
fix a ${\Bbb Z}_2$-grading on a Verma module
$\wdM(\lambda)$ in such a way that a highest weight vector becomes
even. As a $\fg_0$-module,  $\wdM=\wdM_0\oplus\wdM_1$;
the module $\wdM(\lambda)_i$ has a filtration
such that the set of factors coincides
with the set $\{M(\lambda-\gamma):\gamma\in \Gamma_i\}$--- 
see~\cite{mu}, 3.2.

It is easy to check that for any $w\in W,\gamma\in\Gamma$
\begin{equation}
\label{w*}
w.\lambda-w_*\gamma+\rho_0=w(\lambda-\gamma+\rho_0).
\end{equation}
Therefore the $\fg_0$-central characters of
$M(w.\lambda-w_*\gamma)$ and $M(\lambda-\gamma)$
coincide. Thus the multiset of $\fg_0$-central characters 
of $\{M(w.\lambda-\gamma):\gamma\in \Gamma\}$ does not
depend on the choice of $w\in W$. Recall that
the set 
$$W(\wdchi):=\{\mu\in\fh^*|\wdchi \wdM(\mu)=0\}$$
forms a single $W.$-orbit.

\subsubsection{}
\label{bothlie}
Say that $\mu,\mu'\in \sum_i{\Bbb Z}\sigma_i$ have
the same parity if $(\mu-\mu')\in {\Bbb Z}\Delta_0^+$.
Suppose that $\gamma,\gamma'\in\Gamma$ both lie 
either in $\Gamma_0$ or in $\Gamma_1$. This is
equivalent to the condition that the elements
$\gamma,\gamma'$ have the same parity.
For any $w\in W$ one has $w_*\gamma-w_*\gamma'=w(\gamma-\gamma')$
and so the elements $w_*\gamma,w_*\gamma'$ also have the same parity.
Hence both $w_*\gamma,w_*\gamma'$ lie 
either in $\Gamma_0$ or in $\Gamma_1$. 
This has the following consequence:  
the multiset of $\fg_0$-central characters 
of $\{M(w.\lambda-\gamma):\gamma\in \Gamma_0\}$
coincides either with the multiset of $\fg_0$-central characters 
of $\{M(\lambda-\gamma):\gamma\in \Gamma_0\}$
or with the multiset of $\fg_0$-central characters 
of $\{M(\lambda-\gamma):\gamma\in \Gamma_1\}$.
In particular, if $\chi\in\Max\cZ(\fg_0)$ is such
that for some $\lambda\in W(\wdchi)$ one has
$\wdM(\lambda)_{\chi}=M(\lambda-\gamma_0)\oplus M(\lambda-\gamma_1)$
for certain $\gamma_i\in\Gamma_i (i=0,1)$ then
for any  $\lambda'\in W(\wdchi)$ the $\fg_0$-modules
$\wdM(\lambda')_{\chi;i}\,$ are Verma modules.

\subsubsection{}
\label{cnstr}
Define a lexicographic order on ${\Bbb C}$ by setting
$c_1>c_2$ if $\re c_1>\re c_2$ or  $\re c_1=\re c_2$ and $\im c_1>\im c_2$.

The condition~(\ref{condbeta}) implies
the existence of $\lambda\in W(\wdchi)$ such that
$$\lambda+\rho=\sum_1^l k_i\sigma_i,\ \ k_1,\ldots,k_{l-1}>k_l=0.$$
Set $\wdM:=\wdM(\lambda)$ and
$$\chi=\Ann_{\cZ(\fg_0)}M(\lambda).$$ 

One has $\rho_1={1\over 2}\sum_1^l\sigma_i$ and so
$$\lambda+\rho_0=\sum_1^l( k_i+{1\over 2})\sigma_i.$$
One easily sees that 
$$\lambda+\rho_0-\gamma\in W(\lambda+\rho_0)$$
only for $\gamma=\sigma_l, 0$. Therefore $\chi M(\lambda-\gamma)=0$
iff $\gamma=\sigma_l, 0$.
This implies $\wdM_{\chi}=M(\lambda)\oplus M(\lambda-\sigma_l)$.
By~\ref{bothlie}, $\chi$ is a mate for $\wdchi$.

\subsubsection{}
Let us show that $\chi$ is a perfect mate.
Suppose this is not true. 
Then $\wdV(w.\lambda)_{\chi;i}=0$ for some $w\in W$ and $i\in \{0,1\}$.
The equality~(\ref{w*}) implies that for any $y\in W$
$$\wdM(y.\lambda)_{\chi}=M(y.\lambda-y_*0)\oplus M(\lambda-y_*\sigma_l)$$
and so $\wdV(y.\lambda)_{\chi}\ $ is a homomorphic
image of $\ M(y.\lambda-y_*0)\oplus M(\lambda-y_*\sigma_l)$.
The module $\wdV(w.\lambda)$ is a homomorphic
image of the module $\wdM(w.\lambda)$; denote the kernel
of this homomorphism by $\wdN$. 
The module $\wdN$ has finite length and the factors
of its Jordan-G\"older series have form $\wdV(\mu)$
for some $\mu\in W.\lambda$ satisfying $\mu<w.\lambda$.
Since $0=\wdV(w.\lambda)_{\chi;i}=(\wdM(w.\lambda)/\wdN)_{\chi;i}\,$, 
one concludes that the $\fg_0$-module 
$\wdN_{\chi;i}=\wdM(w.\lambda)_{\chi;i}$
has a finite filtration whose factors
are quotients of either 
$M(y.\lambda-y_*0)$ or $M(y.\lambda-y_*\sigma_l)$ for some $y\in W$
satisfying $y.\lambda<w.\lambda$.
Hence the highest weight of $\wdM(w.\lambda)_{\chi;i}\,$
belongs to the set 
$$X:=\{y.\lambda-y_*0,\ y.\lambda-y_*\sigma_l\}_{y\in W\ s.t.\ 
y.\lambda<w.\lambda}.$$
The module $\wdM(w.\lambda)_{\chi;_i}\,$ is isomorphic either to 
$M(w.\lambda-w_*0)$ or to $M(w.\lambda-w_*\sigma_l)$.
Therefore either $w.\lambda-w_*0$ or $w.\lambda-w_*\sigma_l$
belongs  to the set $X$.

If $\,w.\lambda-w_*0=y.\lambda-y_*0\,$ then, by~(\ref{w*}), 
$w(\lambda+\rho_0)=y(\lambda+\rho_0)$, that is 
$w^{-1}y\in\Stab_W(\lambda+\rho_0)$. One has
$$\Stab_W(\lambda+\rho_0)=\Stab_W(\sum_1^l (k_j+{1\over 2})\sigma_j)
\subseteq  \Stab_W(\sum_1^l k_j\sigma_j)=\Stab_W(\lambda+\rho)$$
since  $k_1,\ldots, k_{l-1}>k_l=0$.
Thus $w.\lambda-w_*0=y.\lambda-y_*0$ implies $w.\lambda=y.\lambda$.
Hence 
$(w.\lambda-w_*0)\not\in \{y.\lambda-y_*0\}_{y\in W,y.\lambda<w.\lambda}.$

If $\,w.\lambda-w_*0=y.\lambda-y_*\sigma_l\,$ then,  by~(\ref{w*}),
$w(\lambda+\rho_0)=y(\lambda+\rho_0 -\sigma_l)$ or, equivalently,
$w(\lambda+\rho_0)=ys_{\beta_l}(\lambda+\rho_0)$ where $s_{\beta_l}\in W$
is the reflection with respect to the root $\beta_l$.
Therefore $w^{-1}ys_{\beta_l}\in\Stab_W(\lambda+\rho_0)$.
As we already saw 
$\Stab_W(\lambda+\rho_0)\subseteq  \Stab_W(\lambda+\rho)$  and so
$w^{-1}ys_{\beta_l}\in\Stab_W(\lambda+\rho)$. Then
$w^{-1}y\in\Stab_W(\lambda+\rho)$ since  
$s_{\beta_l}\in\Stab_W(\lambda+\rho)$. Hence
$w.\lambda-w_*0=y.\lambda-y_*\sigma_l$ forces 
$w^{-1}y\in\Stab_W(\lambda+\rho)$. We conclude that
$(w.\lambda-w_*0)\not\in X$.

Similarly, if  $w.\lambda-w_*\sigma_l=y.\lambda-y_*0$ then, by~(\ref{w*}),
$w(\lambda+\rho_0 -\sigma_l)=y(\lambda+\rho_0)$.
As we have shown above, this implies $w^{-1}y\in\Stab_W(\lambda+\rho)$.
Hence $(w.\lambda-w_*\sigma_l)
\not\in \{y.\lambda-y_*0\}_{y\in W,y.\lambda<w.\lambda}.$

Finally, if $w.\lambda-w_*\sigma_l=y.\lambda-y_*\sigma_l$ 
then, by~(\ref{w*}),
$w^{-1}y\in\Stab_W(\lambda+\rho_0 -\sigma_l)$. One can easily deduce
from the equality
$$\lambda+\rho_0-\sigma_l=
\sum_1^{l-1} (k_j+{1\over 2})\sigma_j-{1\over 2}\sigma_l$$
that 
$\Stab_W(\lambda+\rho_0 -\sigma_l)\subseteq  \Stab_W(\lambda+\rho)$.
Thus $(w.\lambda-w_*\sigma_l) \not\in X$ as required.

Hence $\{w.\lambda-w_*0;w.\lambda-w_*\sigma_l\}\cap X=\emptyset$.
This proves that $\wdV(w.\lambda)_{\chi;i}\not=0$ for $i=0,1$.

\subsubsection{}
\begin{cor}{corsemiperfect}
The ideal $\chi\in\Max\cZ(\fg_0)$ described in~\ref{cnstr}
is a perfect mate for $\wdchi$.
\end{cor}

\subsection{}
\label{leftright}
Suppose that $\chi\in\Max\cZ(\fg_0)$ is a perfect mate for $\wdchi$.

The ideal $\Ug\wdchi$ is equal to the annihilator of a Verma module $\wdM$
which has a finite support in $\cZ(\fg_0)$. Therefore
any $\fg$-module annihilated by $\Ug\wdchi$
has a finite support in $\cZ(\fg_0)$. Arguing as in~\cite{g}, 8.3.2,
8.3.3, one can deduce from the definition
of perfect mate that for a graded $\fg$-module $\wdN$ 
annihilated by $\Ug\wdchi$ one has
$\wdN=\Ug\wdN_{\chi;0}$. In particular, $\wdN_{\chi;0}\not=0$ if
$\wdN\not=0$.

For a graded $\fg$-bimodule $L$ set
$$_{\chi}L_{\chi;0}:=\{f\in L_0|\ \chi^r f=f\chi^r=0, \ \forall r>>0\}.$$
Arguing as in~\cite{g}, 8.4 one concludes $_{\chi}L_{\chi;0}\not=0$
provided $L\not=0$ and $\wdchi L=L\wdchi=0$.

\subsection{}
\begin{prop}{propAB}
Let $\wdM$ be a Verma module annihilated by $\wdchi$
and let $A$ be the image of $\Ug$ in $F(\wdM,\wdM)$
under the natural map.
Then 
$$F(\wdM,\wdM)=A\oplus A\theta$$
where $\theta\in F(\wdM,\wdM)$ is given by
$$\theta(v)=(-1)^{d(v)}v.$$ 
\end{prop}
\begin{pf}
Suppose that $\wdM$ is simple. By~\cite{g},
11.1.5, for any simple $\fg$-module $\wdV$ one has 
$\dim\Hom_{\fg}(\wdV,F(\wdM,\wdM))=\dim\Hom_{\fg}(\wdV,H)$
where $H$ is a harmonic space ---see~\ref{annosp12l}. 
The complete reducibility of both $H$ and $F(\wdM,\wdM)$ 
forces $H\cong F(\wdM,\wdM)$. The element $\theta$ is
$\ad\fg_0$-invariant and so $A\cong A\theta$
as $\ad\fg_0$-modules.
In the light of~(\ref{12}), as $\ad\fg_0$-modules
$H\cong A\oplus A$. Since the multiplicity
of any finite dimensional $\fg_0$-module in $H$ is finite,
in order to prove the equality $F(\wdM,\wdM)=A\oplus A\theta$
it is enough to show that $F(\wdM,\wdM)=A+A\theta$.
Since $\theta u=(-1)^{d(u)}u\theta$ for any $u\in\Ug$,
$(A+A\theta)$ is a graded $\Ug$-subbimodule of $F(\wdM,\wdM)$.
Let $\chi\in \cZ(\fg_0)$ be a perfect mate for $\wdchi$.
Then both $M:=\wdM_{\chi;0}$ and $M':=\wdM_{\chi;1}$
are Verma $\fg_0$-modules. The restriction of 
endomorphisms of $\wdM$ to $M\oplus M'$ induces an isomorphism
$$_{\chi}F(\wdM,\wdM)_{\chi;0}\iso F(M,M)\oplus F(M',M').$$
Recall that the natural maps $\cU/(\cU\chi)\to F(M,M)$  
and $\cU/(\cU\chi)\to F(M',M')$ are bijective;
identify $_{\chi}F(\wdM,\wdM)_{\chi;0}$ with  
$\cU/(\cU\chi)\oplus  \cU/(\cU\chi)$ through these maps.
Write
$$\Ann_{\cZ(\fg_0)}\wdM=\chi\prod_1^s \chi_j^{r_j}$$
where $\supp_{\cZ(\fg_0)}=\{\chi,\chi_1,\ldots,\chi_s\}$
and $r_1,\ldots, r_s\in {\Bbb N}^+$. Take $a\in\prod_1^s \chi_j^{r_j}$ 
such that $a=1$ modulo $\chi$. The image of
$\cU a$ in $F(\wdM,\wdM)$ lies in $_{\chi}F(\wdM,\wdM)_{\chi;0}$.
Since for any $u\in\cU$
the element $ua$ acts on both $M$ and $M'$ by the multiplication by $u$,
the image $J$ of $\cU a$  in $_{\chi}F(\wdM,\wdM)_{\chi;0}$
is equal to the diagonal copy of $\cU/(\cU\chi)$ inside
 $\cU/(\cU\chi)\oplus  \cU/(\cU\chi)$:
$$J=\{(u,u)|\ u\in \cU/(\cU\chi)\}.$$
Then $J\theta=\{(u,-u)|\ u\in \cU/(\cU\chi)\}$
and thus 
$$J+J\theta=\cU/(\cU\chi)\oplus  \cU/(\cU\chi)=_{\chi}F(\wdM,\wdM)_{\chi;0}.$$
Hence 
$$_{\chi}F(\wdM,\wdM)_{\chi;0}=_{\chi}(A+A\theta)_{\chi;0}.$$
By~\ref{leftright}, this forces $F(\wdM,\wdM)=A+A\theta$.

We have shown that $F(\wdM,\wdM)=\im\Ug\oplus \theta\im\Ug$ provided 
$\wdM$ is simple.
Take an arbitrary Verma module $\wdM$ satisfying $\wdchi\wdM=0$.
The module $\wdM$ contains a simple Verma submodule
$\wdM'$ since $\Ug$ is a domain. Let $E$ be an $\ad\fg$-submodule
of $\Ug$ such that  the restriction
of the natural map $f':\Ug\to F(\wdM',\wdM')$ to $E$ provides
a bijection $E\to \im f'$. Then $F(\wdM',\wdM')=f'(E)\oplus f'(E)\theta$.
The restriction of the natural map $f:\Ug\to F(\wdM,\wdM)$ to $E$ provides
a bijection $E\to \im f$ since $\ker f=\ker f'=\Ug\wdchi$. 
Suppose that $f(u_1)=f(u_2)\theta$ for some $u_1,u_2\in E$;
then $(u_1-u_2)\wdM_0=0$ and so $f'(u_1)=f'(u_2)\theta$ that implies 
$u_1=u_2=0$ since $f'(E)\cap f'(E)\theta=0$.
Therefore $f(E)\cap f(E)\theta=0$. Thus
$F(\wdM,\wdM)$ contains $f(E)\oplus f(E)\theta$.  A Joseph's reasoning 
based on the use of GK-dimension, shows that the map
$F(\wdM',\wdM')\to F(\wdM',\wdM)$ is bijective and the map
$F(\wdM,\wdM)\to F(\wdM',\wdM)$ is injective (both maps
are induced by the embedding $\wdM'$ to $\wdM$)--- 
see~\cite{jbook}, 8.3.9
or~\cite{g}, 9.2. In particular, $F(\wdM,\wdM)$
is isomorphic to an $\ad\fg$-submodule of $F(\wdM',\wdM')$.
On the other hand, $F(\wdM,\wdM)$
contains $f(E)\oplus f(E)\theta$ which is isomorphic, as an $\ad\fg_0$-module,
to $F(\wdM',\wdM')$. Since the multiplicity of
each finite dimensional $\fg_0$-module in $F(\wdM',\wdM')$ is finite,
one concludes that $F(\wdM,\wdM)=f(E)\oplus f(E)\theta=
\im f\oplus \theta\im f$ as required.
\end{pf}



\begin{thebibliography}{MMMMM}

\bibitem[BF]{bf} A.~D.Bell, R.~Farnsteiner, On the theory
of Frobenius extensions and its application to Lie superalgebras,
Trans AMS {\bf 335} (1993) no.1, p.407---424.      



\bibitem[BL]{bl} J.~Bernstein, V.~Lunts, A simple proof
of Kostant's theorem that $\Ug$ is free over its center,
Amer. J. Math., {\bf 118} (1996), no.5, p. 979--987.

\bibitem[BZV] {bzv}  M.~Bershadsky, S.~Zhukov, A.~Vaintrob, 
${\frak {psl}}(n|n)$ sigma
model as a conformal field theory, Nuclear Phys. B, {\bf 559} (1999),
no.1-2, p. 205--234.

\bibitem[D]{d} M.~Duflo, 
Construction of primitive ideals in an enveloping
algebra, in: I.~M.~Gelfand, ed.. Publ. of 1971 Summer School in Math.,
Janos Bolyai Math. Soc., Budapest p.77---93.



\bibitem[G1]{ghost} M.~Gorelik, 
On the ghost centre of Lie superalgebras, 
to appear at Ann. Inst. Fourier.


\bibitem[G2]{g} M.~Gorelik, Annihilation Theorem and Separation
Theorem for basic classical Lie superalgebras, preprint MSRI 2000-019,
math.RA/0008143.


\bibitem[GL1]{gl2} M.~Gorelik, E.~Lanzmann,
The minimal primitive spectrum of the enveloping 
algebra of the Lie superalgebra ${\frak {osp}}(1,2l)$, to appear at
Adv. in Math.


\bibitem[J1]{j16} A.~Joseph, Kostant's problem, Goldie rank
and the Gelfand-Kirillov conjecture, Invent. Math. {\bf 56} (1980)
p. 193--204. 


\bibitem[J2]{jbook} A.~Joseph, Quantum groups and their primitive
ideals, Springer, 1995.



\bibitem[K1]{kadv} V.~G.~Kac, Lie superalgebras, Adv. in Math. {\bf 26}
(1977) p.8--96.

\bibitem[K2]{kch1} V.~G.~Kac, Characters of typical representations
of Lie superalgebras, Comm. Alg. {\bf 5} (1977) p.889--997.

\bibitem[K3]{kch} V.~G.~Kac, Representations
of classical Lie superalgebras, Lecture Notes in Math., 
Springer-Verlag, Berlin, {\bf 676} 
(1978) p.597--626.


\bibitem[Ko]{ko} B.~Kostant, Lie group representations on polynomial
rings, Amer. J. Math. {\bf 85} (1963) p.327--404.


\bibitem[M1]{mu} I.~M.~Musson, On the center of the enveloping
algebra of a classical simple Lie superalgebra, J. of Algebra,
{\bf 193} (1997), p.75--101.



\bibitem[P]{p} I.~Penkov, Generic representations
of classical Lie superalgebras and their localization, 
Monatshefte f. Math., {\bf 118} (1994) p.267--313.

\bibitem[PS1]{psI} I.~Penkov, V.~Serganova, Representation
of classical Lie superalgebras of type I, Indag. Mathem., N.S. {\bf 3}
(4) (1992), p.419--466.

\bibitem[PS2]{psg} I.~Penkov, V.~Serganova, Generic irreducible
representations of finite-dimensional Lie superalgebras, International
J. Of Math., Vol. 5, No. {\bf 3} (1994) p.389--419.

\bibitem[S1]{s1} A.~N.~Sergeev, Invariant polynomial functions on 
Lie superalgebras, C.R.Acad. Bulgare Sci.
{\bf 35} (1982), no.5, p. 573--576.

\end{thebibliography}
\end{document}